\documentclass[submission]{FPSAC2020}

\usepackage{paralist}
\renewcommand{\th}{\vartheta}

\newcommand{\Egf}{{\sf E}}
\newcommand{\Fgf}{{\sf F}}

\newcommand{\Kgf}{{\sf G}}
\newcommand{\Ktgf}{{\sf \tilde{G}}}

\theoremstyle{plain}
\newtheorem{thm}{Theorem}

\usepackage{lipsum}

\newtheorem*{lem*}{Lemma}
\newtheorem{Theorem}{Theorem}%[section]
\newtheorem{Lemma}[Theorem]{Lemma}

\newtheorem{Corollary}[Theorem]{Corollary}

\graphicspath{{Figures/}}

\title{Counting lattice walks by winding angle}

\author{Andrew Elvey Price\thanks{\href{mailto:andrew.elvey@univ-tours.fr}{This project was supported by the European Research Council (ERC) in the European Union's
Horizon 2020 research and innovation programme, under the Grant Agreement No. 759702.}}\addressmark{1}}

\address{\addressmark{1}LaBRI, Universit\'e de Bordeaux, France}

\received{\today}

\abstract{%We address the problem of counting walks by winding angle around the origin on four different lattices including the square lattice and the triangular lattice. Our method uses a new decomposition of each lattice, which allows us to write functional equations characterising a generating function of walks counted by length, endpoint and winding angle. We then solve these functional equations in terms of Jacobi theta functions.
We address the problem of counting walks by winding angle on the Kreweras lattice, an oriented version of the triangular lattice. Our method uses a new decomposition of the lattice, which allows us to write functional equations characterising a generating function of walks counted by length, endpoint and winding angle. We then solve these functional equations in terms of Jacobi theta functions.
By using this result in conjunction with the reflection principle, we count walks confined to a cone of opening angle any multiple of $\frac{\pi}{3}$, allowing us to extract asymptotic and algebraic information for these walks. %We are also able to extract the limiting distribution of winding angles of walks on each of these lattices, which, as expected coincides with the probability... %One new result in this direction is the enumeration of tandem walks confined to the three quarter plane.
 Our method and results extend analogously to three other lattices, including the square lattice and triangular lattice. On the square lattice, most of our results were derived by Timothy Budd in 2017, so the current work can be seen as an extension of Budd's results to the three other lattices that we consider. Budd's method of deducing these results was very different, as it was based on an explicit eigenvalue decomposition of certain matrices counting paths in the lattice.}


\begin{document}

\maketitle

\section{Introduction}
\begin{figure}[ht]
\centering
   \includegraphics[scale=0.6]{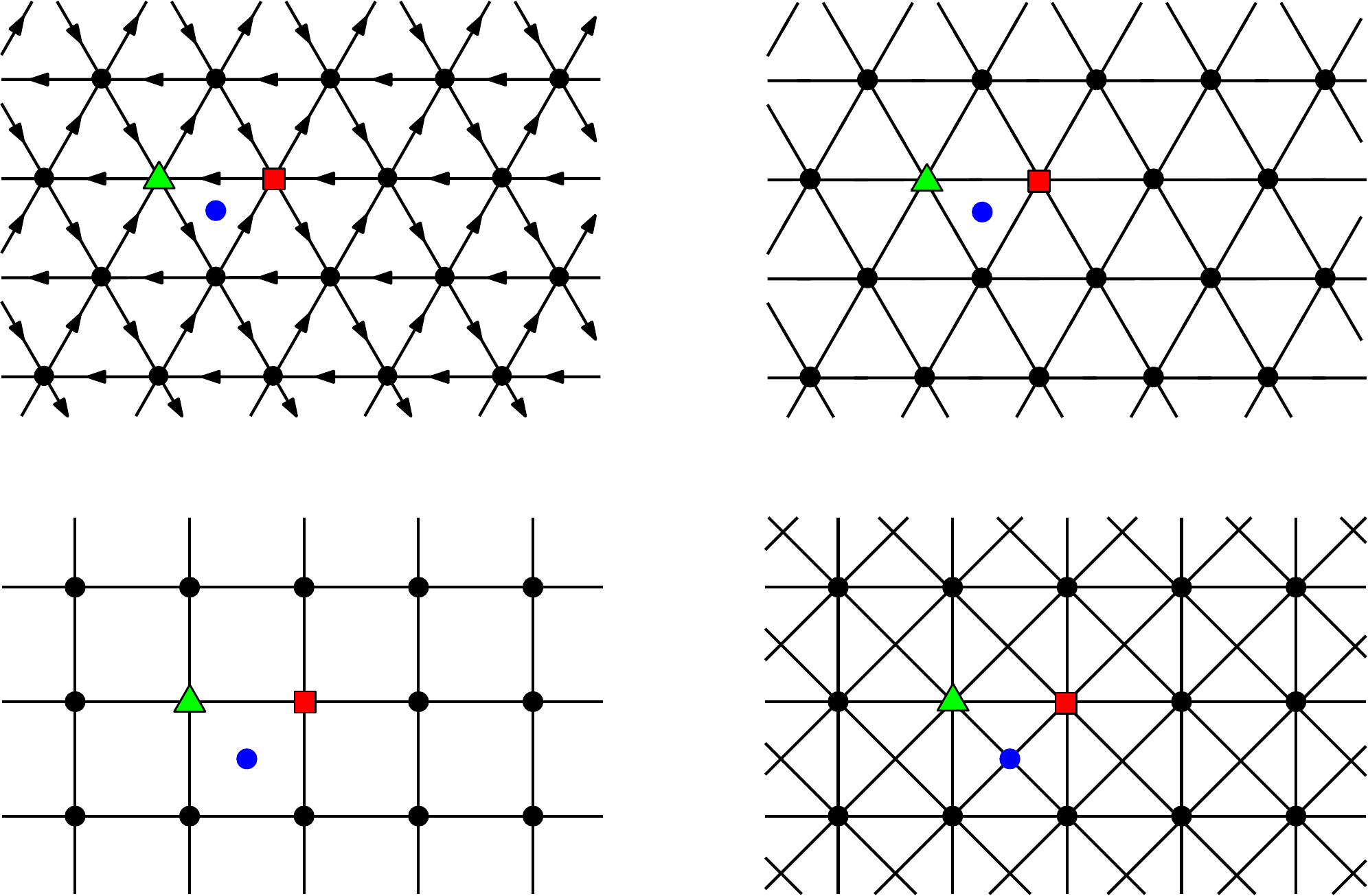}
   \caption{The four lattices on which we count walks by winding number: The Kreweras lattice, the triangular lattice, the square lattice and the king lattice.  Walks may travel in either direction along unoriented edges. In the cell-centred cases the origin is at the blue circle, while in the vertex-centred cases the origin is at the green triangle. In each case, we count walks starting at the red square by winding angle around the origin.}
   \label{fig:lattices}
\end{figure}

We study walks by winding number around the origin $0$ on four different lattices shown in Figure \ref{fig:lattices}: The triangular lattice, the square lattice, the king lattice and the Kreweras lattice. On the square lattice, our results coincide with those derived by Timothy Budd \cite{budd2017winding} using a very different method involving an explicit eigenvalue decomposition of certain matrices.

 Each lattice can be positioned either so that $0$ is at the centre of a cell of the lattice, or so that $0$ is one of the vertices of the lattice. In the former case we call the lattice {\em cell-centred} while in the latter case we call it {\em vertex-centred} (see Figure \ref{fig:lattices}). In each case the walks are forbidden from passing through the origin (this is relevant for the cell-centred king lattice and all vertex-centred lattices). Given a walk $w$, we define the winding angle of $w$ as follows: let $x$ be a variable point that moves continuously along the path $w$, and let $v_{x}=\frac{x}{|x|}$ be a variable unit vector pointing towards $x$. The {\em winding angle} of $w$ is the total anticlockwise angle that $v_{x}$ spins around $0$. If $w$ starts and finishes at the same point, its winding angle $\theta$ will necessarily be a multiple of $2\pi$.% and in this case the {\em winding number} is defined as $\frac{\theta}{2\pi}$.

  For each lattice in Figure \ref{fig:lattices}, in both the vertex-centred and cell-centred cases, we have derived an exact expression counting walks starting at a specific point near the origin by length, endpoint and winding angle. We present here only the results on the Kreweras lattice. %Some of these expressions are quite complicated, so we focus on walks starting {\em and ending} near the origin on the Kreweras lattice.   
%Counting paths by both their endpoint and winding number around $v$ is equivalent to counting paths on the covering space $\Gamma$ of $\mathbb{C}\setminus\{0\}$, shown in Figure \ref{fig:spiral} for the triangular lattice.
All of our results are in terms of the power series $T_{0}(u,q),T_{1}(u,q),\ldots$ defined by
\[T_{k}(u,q)=\sum_{n=0}^{\infty}(-1)^{n}(2n+1)^k q^{n(n+1)/2}(u^{n+1}-(-1)^k u^{-n}).\]
These are related to the Jacobi theta function $\th(z,\tau)\equiv\th_{11}(z,\tau)$ %\equiv\th_{11}(z,\tau)
%\[\th(z,\tau)=\sum_{n=-\infty}^{\infty}(-1)^ne^{\left(\frac{2n+1}{2}\right)^2i\pi\tau+(2n+1)iz}\]
%\\
%&=e^{\frac{\pi\tau i}{4}}\left(e^{iz}-e^{-iz}\right)\prod_{n=1}^{\infty}\left(1-e^{2\pi\tau ni+2iz}\right)\left(1-e^{2\pi\tau n i-2i z}\right)\left(1-e^{2\pi\tau ni}\right)\end{align*}
by
\begin{equation}\label{thandT}\th^{(k)}(z,\tau)=e^{\frac{(\pi\tau-2z)i}{2}}i^{k}T_{k}(e^{2iz},e^{2i\pi\tau}),\end{equation}
where the derivatives are taken with respect to $z$. So, our results can equivalently be written in terms of this theta function.

\begin{figure}[ht!]
\centering
   \includegraphics[scale=0.7]{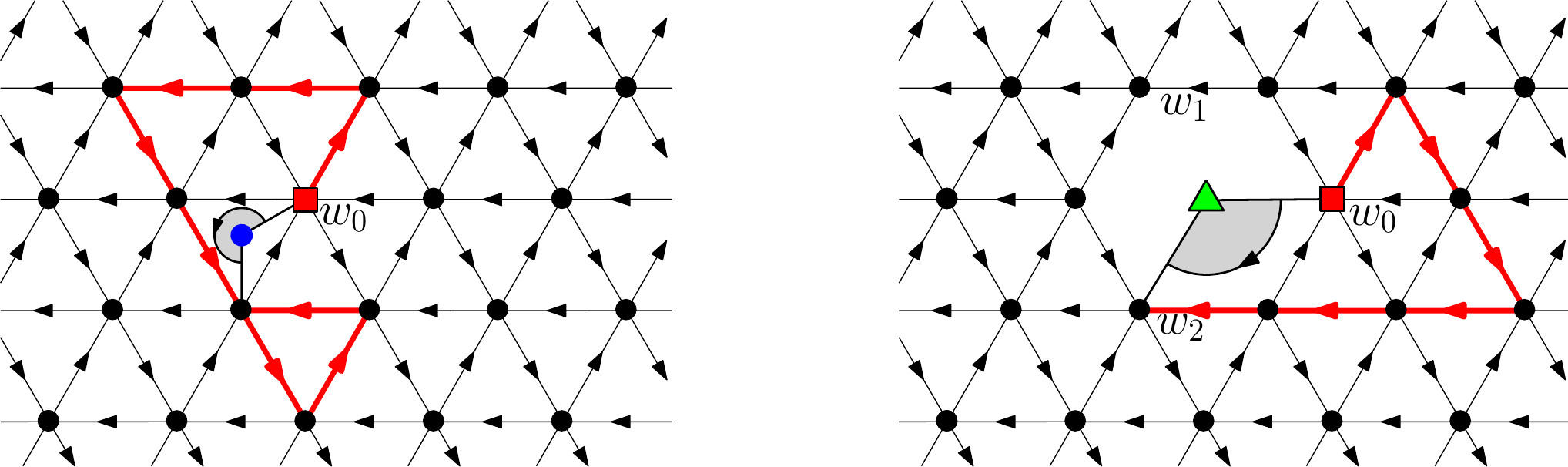}
   \caption{Left: A path with winding angle $\frac{4\pi}{3}$ on the cell-centred Kreweras lattice. Right: a path with winding angle $-\frac{2\pi}{3}$ on the vertex-centred Kreweras lattice.}
   \label{fig:krewerlattice}
\end{figure}

%For our first theorem, let the Kreweras lattice be positioned so that the origin is at the centre of a vertical triangular cell $T$, and let $w_{0}$  be a vertex of this cell, as shown on the left of Figure \ref{fig:krewerlattice}.
For our first theorem we enumerate, by length and winding angle, walks starting and ending near the origin on the cell-centred Kreweras lattice shown on the left of Figure \ref{fig:krewerlattice}. \linebreak
\begin{thm}\label{excursionthm}
Let $q(t)\equiv q=t^3+15t^6+279t^9+\cdots$ be the unique series with constant term 0 satisfying
\[t=q^{1/3}\frac{T_{1}(1,q^3)}{4T_{0}(q,q^3)+6T_{1}(q,q^3)}.\]
The generating function $E(t,s)$ of walks in the cell-centred Kreweras lattice starting at $w_{0}$ and ending at one of the vertices closest to the origin $0$ in which each walk $w$ of length $n$ and winding angle $\frac{2\pi k}{3}$ contributes $t^{n}s^{k}$ is given by
%\[E(t,s)=\frac{s^{1/3}}{1-s}\left(s^{1/3}-q^{-1/3}\frac{T_{1}(q^2,q^3)}{T_{1}(1,q^3)}+q^{-1/3}\frac{T_{0}(q,q^3)T_{1}(s^{1/3}q^{-2/3},q)}{T_{1}(1,q^3)T_{0}(s^{1/3}q^{-2/3},q)}\right)\]
\[E(t,s)=\frac{s}{(1-s^3)t}\left(s-q^{-1/3}\frac{T_{1}(q^2,q^3)}{T_{1}(1,q^3)}-q^{-1/3}\frac{T_{0}(q,q^3)T_{1}(sq^{-2/3},q)}{T_{1}(1,q^3)T_{0}(sq^{-2/3},q)}\right).\]
\end{thm}
The first few terms of the series $E(t,s)$ are given by
\[
E(t,s)=\frac{1}{t}\left(q^{\frac{1}{3}}+sq^{\frac{2}{3}}+\left(s^{2}+\frac{1}{s}\right)q+\cdots\right)=
1+st+\left(s^{2}+\frac{1}{s}\right)t^2+(5+s^{3})t^3+O(t^4).
\]
Note for example that the terms $(5+s^3)t^3$ count the six excursions in the Kreweras lattices of length $3$, of which $5$ have winding angle $0$ and one has winding angle $2\pi$. 
%From the expression above, we can extract a simple summation formula for each $s$-coefficient $E_{c}(t)=[s^c]E(t,s)$:
%\[E_{c}(t)=\]
%For $c\in\mathbb{Z}$, the generating function $E_{c}(t)$ counts Kreweras walks starting at $w_{0}$, ending at a vertex adjacent to $v$ and having winding angle $\frac{2\pi c}{3}$.

For our second theorem we enumerate, by length and winding angle, certain walks on the vertex-centred Kreweras lattice shown on the right of Figure \ref{fig:krewerlattice}.

%For our second theorem, we enumerate walks on the vertex-centred Kreweras lattice so that the origin is one of the vertices of the lattice, and we let $w_{0}=1$ be the vertex to the right of the origin. Let $\tilde{E}(t,s)$ be the generating function for walks starting at $w_{0}$ and ending at one of the cubic roots of unity in which each walk $w$ of length $n$ and winding angle $\theta$ contributes $t^{n}s^{\frac{3\theta}{2\pi}}$.
\begin{thm}\label{mvertthm}
Let $w_{0},w_{1},w_{2}$ be the three vertices on the vertex-centred Kreweras lattice from which there is an edge pointing to $0$. 
Let $q(t)\equiv q=t^3+15t^6+279t^9+\cdots$ be the same series as in Theorem \ref{excursionthm}.
Then the generating function $\tilde{E}(t,s)$ of walks in the vertex-centred Kreweras lattice starting at $w_{0}$ and ending at one of $w_{0}$, $w_{1}$, $w_{2}$ in which each walk $w$ of length $n$ and winding angle $\frac{2\pi k}{3}$ contributes $t^{n}s^{k}$ is given by
%\[E(t,s)=\frac{s^{1/3}}{1-s}\left(s^{1/3}-q^{-1/3}\frac{T_{1}(q^2,q^3)}{T_{1}(1,q^3)}+q^{-1/3}\frac{T_{0}(q,q^3)T_{1}(s^{1/3}q^{-2/3},q)}{T_{1}(1,q^3)T_{0}(s^{1/3}q^{-2/3},q)}\right)\]
\[\tilde{E}(t,s)=\frac{s q^{-2/3}}{t(1+s+s^2)}\frac{T_{0}(q,q^3)^2}{T_{1}(1,q^3)^2}\left(\frac{T_{1}(q,q^{3})^2}{T_{0}(q,q^{3})^2}-\frac{T_{2}(q,q^{3})}{T_{0}(q,q^{3})}-\frac{T_{2}(s,q)}{2T_{0}(s,q)}+\frac{T_{3}(1,q)}{6T_{1}(1,q)}+\frac{T_{3}(1,q^3)}{3T_{1}(1,q^3)}\right).\]
\end{thm}
Analysing the solutions for $E(t,s)$ and $\tilde{E}(t,s)$ leads to the following corollaries regarding the algebraic and asymptotic nature of these series. Analogous results to these corollaries on square lattice are, respectively, Corollary 20 and Lemma 22 in \cite{budd2017winding}.
\begin{Corollary}\label{cor:albegraicity}
If $s$ is a root of unity which is not a cubic root of unity, then the generating functions $E(t,s)$ and $\tilde{E}(t,s)$ are algebraic in $t$.
\end{Corollary}

\begin{Corollary}\label{cor:asymptotics}
For $s=e^{i\alpha}$, with $\alpha\in(0,\pi)\setminus\left\{\frac{2\pi}{3}\right\}$, % the coefficients $c_{n}=[t^{n}]E(t,s)$ behave like
%\[c_{n}\sim\frac{3^{\frac{3}{2}-\frac{3\alpha}{\pi}}e^{\alpha i-\frac{\pi i}{6}}}{(1-e^{3\alpha i})\Gamma(-\frac{3\alpha}{2\pi})}n^{-\frac{3\alpha}{2\pi}-1}3^{n}\]
%for large $n$, while
 the coefficients $v_{n}=[t^{n}]\tilde{E}(t,s)$ satisfy
\[v_{n}\sim-\frac{3^{5-\frac{3\alpha}{\pi}}e^{\alpha i}\alpha}{2\pi(1+e^{\alpha i}+e^{2\alpha i})\Gamma(-\frac{3\alpha}{2\pi})}n^{-\frac{3\alpha}{2\pi}-1}3^{n},\]
for large $n$ satisfying $3|n$. When $3\nmid n$, we have $v_{n}=0$.
\end{Corollary}

%We delay the statement of our most general results to sections \ref{sec:Krewerwinding} and \ref{sec:otherlattices}, after we have introduced the generating function that we will use to count these walks in general.

In Section \ref{sec:Krewerwinding} we outline the proof of these theorems, and also give a more general result counting all walks by endpoint and winding angle. In Section \ref{sec:walksonacone} we use these results to count walks on the Kreweras lattice confined to a cone of opening angle $\frac{\pi N}{3}$ for any $N\in\mathbb{N}$. We then compute the asymptotic number of these walks as the length approaches infinity, which we find to be consistent with the probabilistic results of \cite{denisov2015random}. We then show that the associated generating function is algebraic if and only if $3\nmid N$. Finally, in Section \ref{sec:otherlattices} we briefly discuss the other three lattices.

\section{Walks by winding number on the Kreweras lattice}
\label{sec:Krewerwinding}

\subsection{Decomposition and functional equations}
\begin{figure}[ht!]
\centering
\begin{picture}(350,180)
   \put(0,0){\includegraphics[scale=0.3]{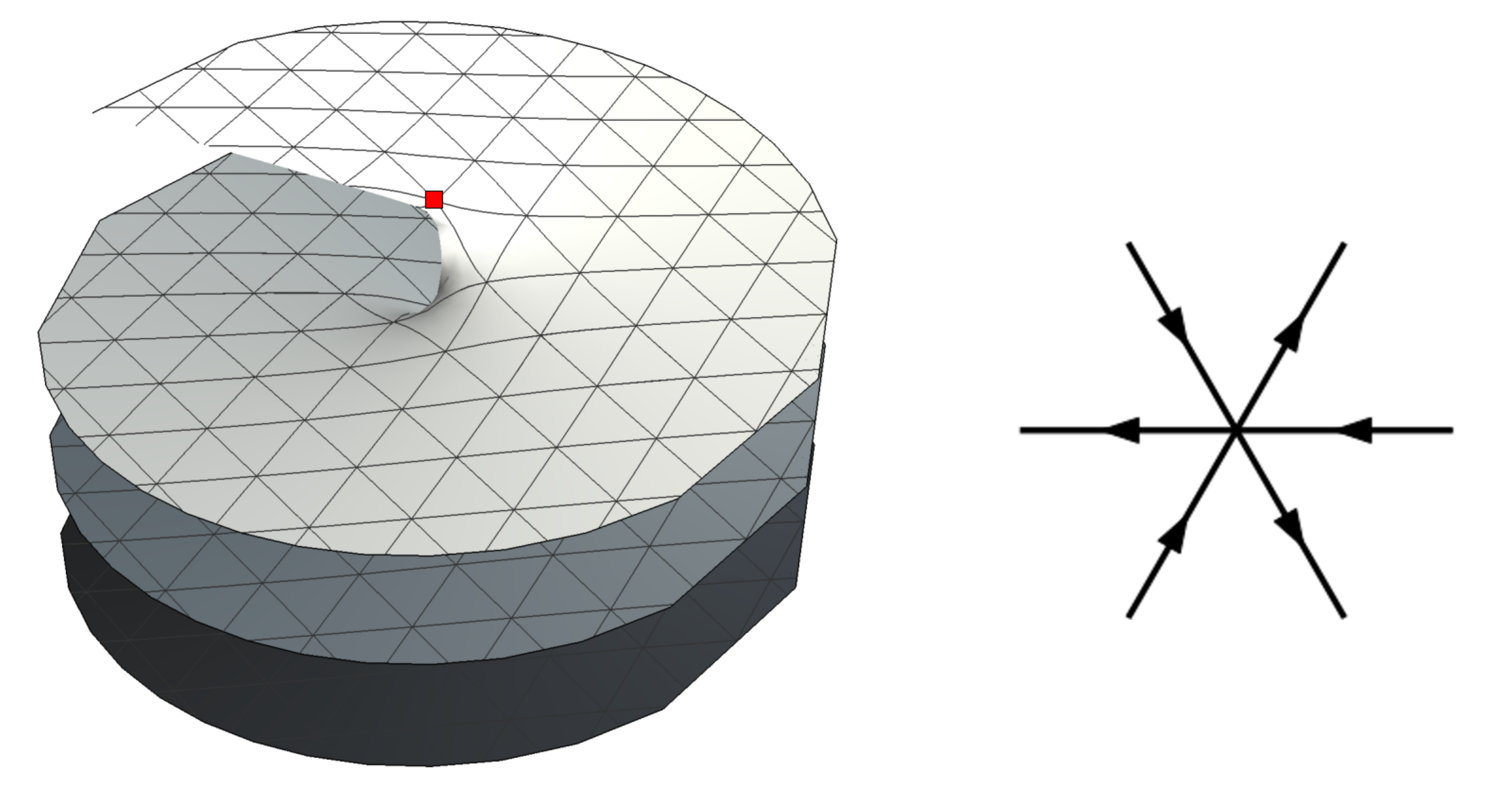}}
   \put(10,10){\large $\Gamma$}
   \put(102,138){\textcolor{red}{\large $w_{0}$}}
   %\put(350,200){$\Gamma$}
\end{picture}
   \caption{Left: The covering space $\Gamma$ of $\mathbb{C}\setminus\{0\}$. Right: The allowed steps for Kreweras walks. We enumerate walks on $\Gamma$, starting at the red square and following these directions.}
   \label{fig:spiral}
\end{figure}

\begin{figure}[ht]
\centering
\includegraphics[scale=0.7]{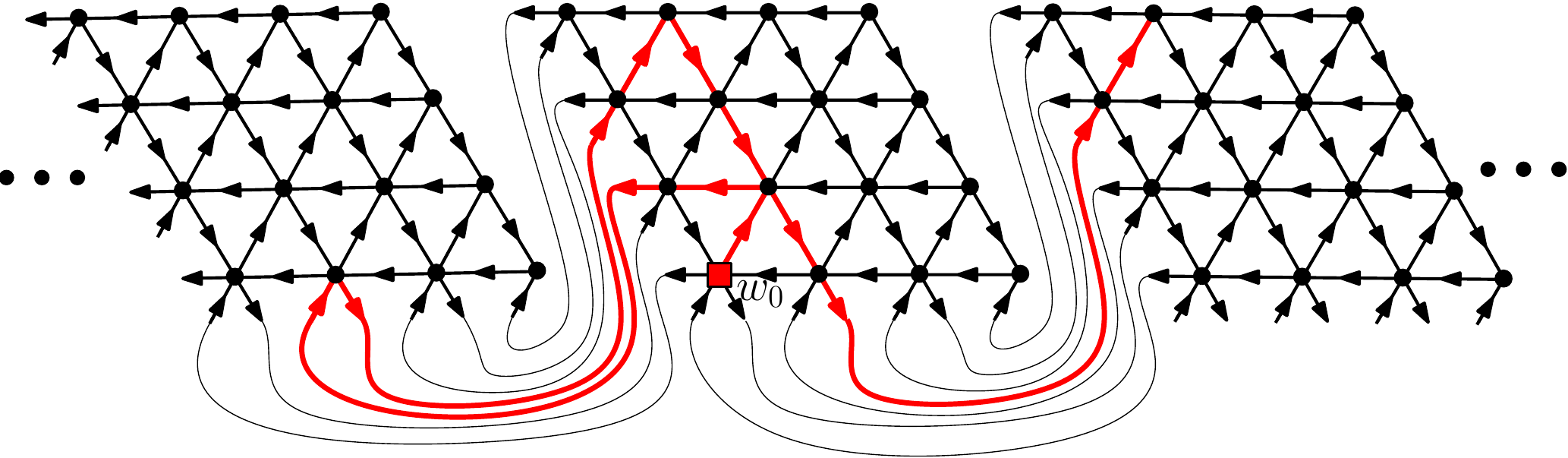}
   \caption{Decomposition of the covering space $\Gamma$ into an infinite sequence of wedges.}
   \label{fig:wedgewalk}
\end{figure}

\begin{figure}[ht]
\centering
\includegraphics[scale=0.7]{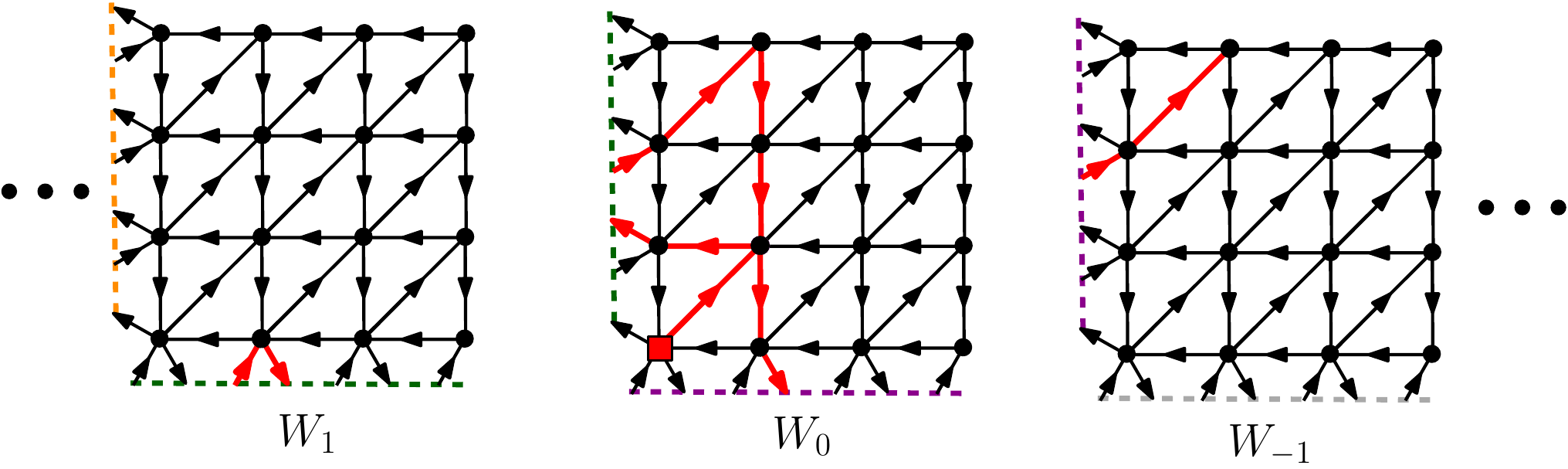}
   \caption{An example of a walk on $\Gamma$ using Kreweras steps. This walk contributes $t^{10}xy^3s^{-1}$ to $\Kgf(t,s,x,y)$.}
   \label{fig:krewerwalk}
\end{figure}

Counting walks in the plane by endpoint and winding number on the cell-centred Kreweras lattice is equivalent to counting walks in the covering space $\Gamma$ of $\mathbb{C}\setminus\{0\}$ by endpoint (see Figure \ref{fig:spiral}). For each point $v\in\Gamma$ and each $n\in\mathbb{N}_{0}$, let $p_{n,v}$ be the number of paths of length $n$ in $\Gamma$, using Kreweras steps, starting at $w_{0}$ and ending at $v$. In this section we study a generating function containing every value of $p_{n,v}$. 
In order to characterise the points $v\in\Gamma$, we partition $\Gamma$ into an infinite sequence of wedges, as shown in Figure \ref{fig:wedgewalk}, then we transform these wedges into quarter-planes $\{W_{j}\}_{j\in\mathbb{Z}}$ as in Figure \ref{fig:krewerwalk}. Each possible endpoint $v$ then corresponds to a triple $(a,b,k)$, where $v\in W_{k}$ and $(a,b)$ are the coordinates of $v$ in $W_{k}$. This way we associate each point $v\in\Gamma$ with a unique monomial $f_{v}(s,x,y)=s^{k} x^a y^b $. The most general series $\Kgf(t,s,x,y)\in\mathbb{Z}\left[s,s^{-1},x,y\right][[t]]$ that we study is
\[\Kgf(t,s,x,y)=\sum_{n=0}^{\infty}t^{n}\sum_{v\in\Gamma}p_{n,v}f_{v}(s,x,y).\]
This is related to the series $\Egf(t,s)$ from Theorem \ref{excursionthm} by $\Egf(t,s)=\Kgf(t,s,0,0).$

By considering all possible final steps of a path in $\Gamma$, we derive the following functional equation, which characterises $\Kgf(x,y)\equiv \Kgf(t,s,x,y)\in\mathbb{Z}\left[s,x,y\right][[t]]$:
\begin{align}\label{eq:krewer}\begin{split}\Kgf(x,y)&=1+txy\Kgf(x,y)+\frac{t}{x}(\Kgf(x,y)-\Kgf(0,y))+\frac{t}{y}(\Kgf(x,y)-\Kgf(x,0))\\
&~+ts\Kgf(0,x)+ts^{-1}y\Kgf(y,0).\end{split}\end{align}
The last two terms in this equation come from the steps from one wedge to an adjacent wedge, so removing these terms yields the well known functional equation for the generating function $\Fgf(t,x,y)$ of paths confined to a quarter plane \cite{bousquet2010walks} with the step set $\{(1,1),(0,-1),(-1,0)\}$ (sometimes called the {\em Kreweras} step set).

\subsection{Solution in terms of theta functions}
In order to solve the functional equation \eqref{eq:krewer}, we fix $t,s\in\mathbb{C}$ and think of $\Kgf(x,y)$ as an analytic function of $x$ and $y$. More precisely, we assume $|t|<1/3$ and $|s|=1$, as this ensures that the series defining $\Kgf(x,y)\equiv \Kgf(t,s,x,y)$ converges absolutely for $|x|,|y|<1$, using the fact that the total number of paths of length $n$ in $\Gamma$ is $3^n$.

The next step is to analyse the kernel
\[K(x,y)=1-txy-\frac{t}{x}-\frac{t}{y},\]
which is the factor of $\Kgf(x,y)$ in \eqref{eq:krewer}.
The kernel here is exactly the same as in the case for paths confined to a quarter plane with step set $\{(1,1),(0,-1),(-1,0)\}$. Moreover, Kurkova and Raschel  \cite{kurkova2012functions} showed that the kernel equation $K(x,y)=0$ can be parameterised in terms of the Weierstrass elliptic function $\wp$ for almost any step set $S\subset\{-1,0,1\}^2$. In this case, as well as the other step sets that we consider, the parameterisation can be written more simply in terms of the Jacobi theta function
\[\th(z,\tau)=\sum_{n=-\infty}^{\infty}(-1)^ne^{\left(\frac{2n+1}{2}\right)^2i\pi\tau+(2n+1)iz}\]
as follows:
\begin{Lemma}
Let $\tau\in i\mathbb{R}_{>0}$ satisfy
\[t=e^{-\frac{\pi\tau i}{3}}\frac{\th'(0,3\tau)}{4i\th(\pi\tau,3\tau)+6\th'(\pi\tau,3\tau)}\] 
and let $X(z)$ and $Y(z)$ be defined by
\[X(z)=e^{-\frac{4\pi\tau i}{3}}\frac{\th(z,3\tau)\th\left(z-\pi\tau,3\tau\right)}{\th\left(z+\pi\tau,3\tau\right)\th\left(z-2\pi\tau,3\tau\right)}~~~\text{and}~~~Y(z)=X(z+\pi\tau).\]
Then $K(X(z),Y(z))=0$ for all $z\in\mathbb{C}$.
\end{Lemma}
To prove this lemma it suffices to check that the function $\Omega(z)=K(X(z),Y(z))$ is doubly periodic, has no poles and that $\Omega(0)=0$.
%\begin{itemize}
%\item $\Omega(z)$ is elliptic, with $\pi$ and $3\pi\tau$ as periods, 
%\item $\Omega(z)$ has no poles,
%\item $\Omega(0)=0$.
%\end{itemize}
Then the result $\Omega(z)=0$ follows from Liouville's theorem, which states that the only holomorphic, doubly periodic functions are constant.% The first of these properties is a direct consequence of the following pseudo-periodicity relations of $\th(z,\tau)$:
%\[\th(z+\pi,\tau)=-\th(z,\tau)~~~\text{and}~~~\th(z+\pi\tau,\tau)=-e^{-2iz-\pi\tau i}\th(z,\tau).\] 

Since $X(0)=Y(0)=0$, we can substitute $x=X(z)$ and $y=Y(z)$ into \eqref{eq:krewer} for $z$ in a neighbourhood of $0$. Writing
\begin{equation}\label{Ldefn}L(z)=st\Kgf(0,X(z))-\frac{t}{Y(z)}\Kgf(X(z),0),\end{equation}
the resulting equation simplifies to
\begin{equation}\label{Leq}1=-L(z)+\frac{L(z+\pi\tau)}{sX(z)}.\end{equation}
Moreover, we can prove that the only poles of $L(z)$ in $0\leq\text{Im}(z)\leq\text{Im}(\pi\tau)$ are simple poles at the points in $\pi\mathbb{Z}$. Together with \eqref{Leq}, this uniquely defines $L(z)$:
\begin{align*}
L(z)&=\frac{1}{1-e^{3i\alpha}}\left(e^{3i\alpha}+\frac{e^{2i\alpha}}{X(z)}+e^{i\alpha}X(z-\pi\tau)\right)\\
&~+\frac{e^{i\alpha+\frac{5i\pi\tau}{3}}\th(\pi\tau,3\tau)\th'(0,\tau)}{(1-e^{3i\alpha})\th(\frac{\alpha}{2}-\frac{2\pi\tau}{3},\tau)\th'(0,3\tau)}\frac{\th(z-2\pi\tau,3\tau)\th(z-\frac{\alpha}{2}+\frac{2\pi\tau}{3},\tau)}{\th(z,\tau)\th(z,3\tau)}.
\end{align*}
We can then extract $\Kgf(X(z),0)$ from \eqref{Ldefn} using the fact that $X(z)=X(\pi\tau-z)$:
\[\Kgf(X(z),0)=\frac{L(z)-L(\pi\tau-z)}{tX(z)(X(z+\pi\tau)-X(z-\pi\tau))}.\]
Hence we have an exact, parametric expression for $\Kgf(x,0)$. We can similarly derive such an expression for $\Kgf(0,y)$, which yields the exact solution for $\Kgf(x,y)$, using \eqref{eq:krewer}. Substituting $z\to0$ yields the following expression for $\Egf(t,s)=\Kgf(0,0)$:
\[E(t,e^{i\alpha})=\frac{e^{i\alpha}}{t(1-e^{3i\alpha})}\left(e^{i\alpha}-e^{\frac{4\pi\tau i}{3}}\frac{\th'(2\pi\tau,3\tau)}{\th'(0,3\tau)}-e^{\frac{\pi\tau i}{3}}\frac{\th(\pi\tau,3\tau)\th'(\frac{\alpha}{2}-\frac{2\pi\tau}{3},\tau)}{\th'(0,3\tau)\th(\frac{\alpha}{2}-\frac{2\pi\tau}{3},\tau)}\right).\]
%\notea{I'm a bit confused about whether the second $-$ above should be $+$.}
Theorem \ref{excursionthm} follows by expanding both $t$ and $E(t,e^{i\alpha})$ as series in $s=e^{i\alpha}$ and $q=e^{2\pi\tau i}$.

On the vertex centred Kreweras lattice we have a slightly different decomposition, which leads to the functional equation
\begin{align}\label{eq:krewerhole}\begin{split}\Ktgf(x,y)&=1+txy\Ktgf(x,y)+\frac{t}{x}(\Ktgf(x,y)-\Ktgf(0,y))+\frac{t}{y}(\Ktgf(x,y)-\Ktgf(x,0))\\
&~+\frac{ts}{x}\left(\Kgf(0,x)-\Ktgf(0,0)\right)+ts^{-1}y^2\Ktgf(y,0).\end{split}\end{align}
Using the same method as above, we solve this functional equation and, in doing so, derive Theorem \ref{mvertthm}.

\subsection{Asymptotic results}
For fixed $\alpha\in(0,\pi)\setminus\{\frac{2\pi}{3}\}$, we determine the asymptotic behaviour of the coefficients $[t^n]\tilde{E}(t,e^{i\alpha})$
for large $n$. The first step is to use the Jacobi identity
\[\th(z,\tau)=i(-i\tau)^{-1/2}e^{-\frac{i}{\pi\tau}z^2}\th\left(\frac{z}{\tau},-\frac{1}{\tau}\right)\]
to write $t$ and $\tilde{E}(t,e^{i\alpha})$ in terms of $\hat{\tau}=-\frac{1}{3\tau}$. For $t$, the expression is
$t=\frac{\th'(0,\hat\tau)}{6\th'(\pi/3,\hat\tau)}$. 
%tE(t,e^{i\alpha})&=\frac{e^{i\alpha}}{1-e^{3i\alpha}}\left(e^{i\alpha}+\frac{\frac{3i\alpha}{\pi}\th(\frac{\pi}{3},\hat{\tau})+\th'(\frac{\pi}{3},\hat{\tau})}{\th'(0,\hat{\tau})}+\frac{3\th(\frac{\pi}{3},\hat{\tau})\th'(\frac{2\pi}{3}+\frac{3\alpha\hat{\tau}}{2},3\hat{\tau})}{\th'(0,\hat{\tau})\th(\frac{2\pi}{3}+\frac{3\alpha\hat{\tau}}{2},3\hat{\tau})}\right)
%however the expression for $\tilde{E}(t,e^{i\alpha})$ is quite complicated.
Then using \eqref{thandT}, both $t$ and $\tilde{E}(t,e^{i\alpha})$ can be expanded as series in $\hat{q}=e^{2\pi i\hat{\tau}}$:
\begin{align*}
t&=\frac{1}{3}-3\hat{q}+18\hat{q}^2+O(\hat{q}^3)\\
t\tilde{E}(t,e^{i\alpha})&=a_{0}+a_{1}\hat{q}-\frac{27\alpha e^{i\alpha}}{2\pi(1+e^{i\alpha}+e^{2i\alpha})}\hat{q}^{\frac{3\alpha}{2\pi}}-\frac{27(2\pi-\alpha)e^{i\alpha}}{2\pi(1+e^{i\alpha}+e^{2i\alpha})}\hat{q}^{3-\frac{3\alpha}{2\pi}}+a_{2}\hat{q}^{\frac{3\alpha}{2\pi}+1}+O(\hat{q}^2),
\end{align*}
where $a_{0}$, $a_{1}$ and $a_{2}$ are explicit functions of $\alpha$. We can prove that the dominant singularity $t=1/3$ corresponds to $\hat{q}=0$, and from the series expansion for $t$ above, we have $\hat{q}\sim \frac{1-3t}{9}$ in the vicinity of this point. The behaviour of the coefficients of $\tilde{E}(t,s)$ is determined by the term
\[-\frac{27\alpha e^{i\alpha}}{2\pi(1+e^{i\alpha}+e^{2i\alpha})}\hat{q}^{\frac{3\alpha}{2\pi}}\sim -\frac{3^{3-\frac{3\alpha}{\pi}}\alpha e^{i\alpha}}{2\pi(1+e^{i\alpha}+e^{2i\alpha})}(1-3t)^{\frac{3\alpha}{2\pi}}.\]
After taking care of the singularities of $\tilde{E}(t,s)$ at $t=\frac{1}{3}e^{\frac{2i\pi}{3}}$ and $t=\frac{1}{3}e^{-\frac{2i\pi}{3}}$, this yields Corollary \ref{cor:asymptotics}.

In the cases $\alpha=0$ and $\alpha=\frac{2\pi}{3}$, the numerator and denominator in our expression for $\tilde{E}(t,e^{i\alpha})$ are both $0$, so we need to take derivatives in $\alpha$ to obtain the correct expression. We then apply the same method as above to determine the asymptotic form of $v_{n}$. We find that $\tilde{E}(t,e^{i\alpha})$ exhibits a logarithmic singularity in these cases.

%\notea{To complete!! (we will get the asymptotic form, which agrees with the probability stuff and we also get the constant)} 

\subsection{Algebraicity results}
Due to the relation to the theta function, the series $\Kgf(t,s,x,y)$ and $\Ktgf(t,s,x,y)$ are {\em differentially algebraic} in each variable, meaning that they satisfy a non-trivial algebraic differential equation with respect to that variable. For certain specialisations, these series take a much simpler form:
\begin{Corollary} Let $N$ be a positive integer and let $s$ satisfy $s^{N}=1$ but $s^3\neq1$. The generating functions $\Kgf(x,y)\equiv\Kgf(t,s,x,y)$ and $\Ktgf(t,s,x,y)$ are algebraic in $t$, $x$ and $y$.\end{Corollary}
The algebraicity of $\Kgf(x,0)$ in $x$ follows from the fact that the functions $X(z)$ and $\Kgf(X(z),0)$ are both elliptic with periods $\pi$ and $6N\pi\tau$, since any two elliptic functions with the same periods are algebraically related. For the algebraicity of $\Kgf(x,0)$ in $t$, we use modular properties of $t$ and $\Kgf(t,s,X(z),0)$ as functions of $\tau$. The proof that $\Kgf(0,y)$ is algebraic in $t$ and $y$ is similar, and together these imply that $\Kgf(t,s,x,y)$ is algebraic in $t$, $x$ and $y$.

In the case $s^3=1$, even $\Egf(t,s)=\Kgf(t,s,0,0)$ and $\tilde{\Egf}(t,s)$ are not algebraic, as they each have a logarithmic singularity. Nonetheless, we can use similar ideas as above to show that $\Egf(t,s)$ and $\tilde{\Egf}(t,s)$ are D-finite in these cases, meaning they each satisfy a non-trivial linear differential equation.

%Since $\Kgf(t,s,x,y)$ is a rational function of $t$, $x$, $y$, $\Kgf(t,s,x,0)$ and $\Kgf(t,s,0,y)$, The functions $X(z)$, $Y(z)$, $L(z)$, $\Kgf(X(z),0)$ and $\Kgf(0,Y(z))$ are all elliptic with periods $\pi$ and $6N\pi\tau$, meaning that every pair of them is algebraically related. In particular $\Kgf(x,0)$ and $\Kgf(0,y)$ are algebraic in $x$ and $y$. It follows that $\Kgf(x,y)$ is algebraic in $x$ and $y$. There is therefore polynomial $P$ whose coefficients are series in $t$ satisfying $P(. To see the algebraicity in $t$, we observe that, 

%\section{Distribution of winding angles}

\section{Walks confined to a cone}
\label{sec:walksonacone}
Walks in cones are the subject of a wealth of literature in both combinatorics \cite{bousquet2010walks,bousquet2016square} and probability theory \cite{fayolle1999random,denisov2015random}. The simplest non-trivial example is that of walks in the quarter-plane. The systematic study of these walks was started by Bousquet-M\'elou and Mishna in 2010 \cite{bousquet2010walks}, where they considered walks with any step set $S\subset\{-1,0,1\}^{2}$, showing that there are 79 non-trivial and combinatorially distinct such step sets. This study is now in some sense complete as differentially algebraic solutions for the generating function are now known for 32 step sets $S$, while for the other $47$ step sets, the generating function is known to be non-differentially algebraic \cite{dreyfus2018nature}. %Perhaps the next two most natural cones to consider are the slit plane \cite{bousquet2002walks} and
In recent years several authors have started to extend this by considering walks with small steps on the three quarter plane \cite{bousquet2016square,raschel2018walks,bousquet_wallner}. Another interesting cone on which simple walks have been counted is the slit plane \cite{bousquet2002walks}, that is the cone with opening angle $2\pi$.

In this section we use our results, in conjunction with the reflection principle, to count Kreweras walks on a cone with opening angle $\frac{N\pi}{3}$ for any positive integer $N$. Analogous results were found by Budd on the square lattice \cite[Section 3.3]{budd2017winding}. One new result of this type is the enumeration of walks confined to the three quarter plane using the step set $\{(0,1),(-1,0),(1,-1)\}$.

\begin{figure}[ht]
\centering
\includegraphics[scale=0.7]{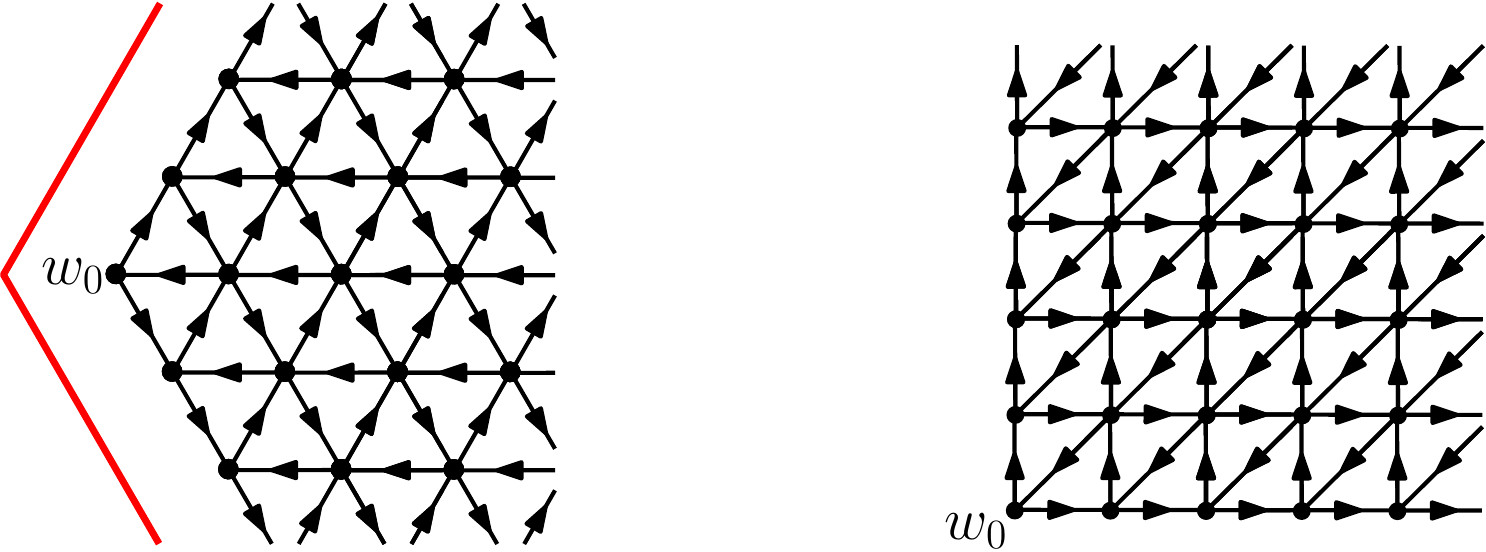}
   \caption{The Kreweras lattice with opening angle $\frac{2\pi}{3}$ is equivalent to the quarter plane with the step set $\{(-1,-1),(0,1),(1,0)\}$. Excursions from $w_{0}$ are counted by $\tilde{E}_{0,-1,1}(t)$.}
   \label{fig:1on3}
\end{figure}

\begin{figure}[ht]
\centering
\includegraphics[scale=0.7]{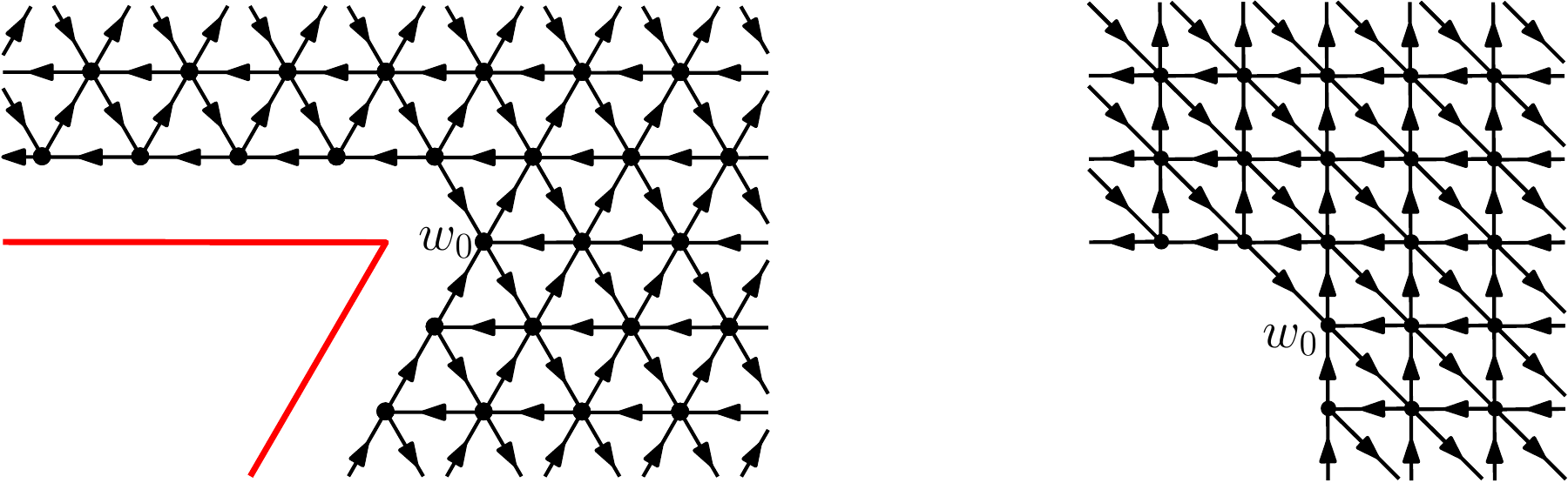}
   \caption{The Kreweras lattice with opening angle $\frac{5\pi}{6}$ is equivalent to the three quarter plane with step set $\{(0,1),(-1,0),(1,-1)\}$. Excursions from $w_{0}$ are counted by $\tilde{E}_{0,-2,3}(t)$.}
   \label{fig:5on6}
\end{figure}

For integers $k_{1}<0<k_{2}$ and $k$ satisfying $k_{1}<2k<k_{2}$, let $\tilde{E}_{k,k_{1},k_{2}}(t)$ be the generating function for Kreweras walks $w$ with the following properties:
\begin{compactitem}
\item $w$ starts at a fixed vertex $w_{0}$, adjacent to $0$,
\item $w$ ends at a vertex adjacent to $0$,
\item $w$ does not pass through $0$,
\item $w$ has winding angle $\frac{2k\pi}{3}$,
\item The winding angle is confined to stay within the interval $(\frac{k_{1}\pi}{3},\frac{k_{2}\pi}{3})$ for all intermediate points along $w$.
\end{compactitem}
Let $\tilde{E}_{k}(t)$ denote the generating function of Kreweras walks satisfying only the first three properties. Then the generating function $\tilde{E}(t,s)$ from Theorem \ref{mvertthm} is given by
\[\tilde{E}(t,s)=\sum_{k=-\infty}^{\infty}s^k\tilde{E}_{k}(t).\]Using the reflection principle we can write $\tilde{E}_{k,k_{1},k_{2}}(t)$ in terms of $\tilde{E}_{k}(t)$:
\[\tilde{E}_{k,k_{1},k_{2}}(t)=\sum_{j=-\infty}^{\infty}\left(\tilde{E}_{k-jk_{1}+jk_{2}}(t)-\tilde{E}_{-k-jk_{1}+(j+1)k_{2}}(t)\right).\]
This can be written as a sum of values of $\tilde{E}(t,s)$ as follows:
\begin{equation}\label{eq:walksinconeeq}\tilde{E}_{k,k_{1},k_{2}}(t)=\frac{1}{k_{2}-k_{1}}\sum_{j=1}^{k_{2}-k_{1}-1}\left(e^{-\frac{2\pi ijk}{k_{2}-k_{1}}}-e^{\frac{2\pi ij(k-k_{1})}{k_{2}-k_{1}}}\right)\tilde{E}\left(t,e^{\frac{2\pi ij}{k_{2}-k_{1}}}\right).\end{equation}
\textbf{Example 1:} By \eqref{eq:walksinconeeq}, the generating function $\tilde{E}_{0,-1,1}(t)$ is given by
\[\tilde{E}_{0,-1,1}(t)=\tilde{E}(t,-1).\]
This counts walks starting and ending at $(1,0)$, which are confined to stay within the wedge with opening angle $120^{\circ}$, shown in Figure $\ref{fig:1on3}$. Equivalently, $E_{0,-1,1}(t)$ counts quarter-plane excursions with the step set $\{(-1,-1),(0,1),(1,0)\}$ (these are sometimes called {\em reverse Kreweras walks}). Note that Corollary \ref{cor:albegraicity} implies that this generating function is algebraic. Indeed this step set is one of only four for which the generating function of excursions on the quarter-plane is algebraic \cite{bousquet2010walks}.
~\newline
\textbf{Example 2:} By \eqref{eq:walksinconeeq}, the generating function $E_{0,-2,3}(t)$ is given by
\[\tilde{E}_{0,-2,3}(t)=\frac{1}{5}\sum_{j=1}^{4}\left(1-e^{\frac{4\pi ij}{5}}\right)\tilde{E}\left(t,e^{\frac{2\pi i}{5}}\right)\]
 counts walks in the $5/6$-plane. Equivalently, this counts walks in the three quarter plane, starting and ending at $(0,-1)$, using the step-set $\{(0,1),(-1,0),(1,-1)\}$ (see Figure \ref{fig:5on6}). Again by Corollary \ref{cor:albegraicity}, this generating function is algebraic in $t$.

In these two examples we saw that $\tilde{E}_{k,k_{1},k_{2}}(t)$ is algebraic in $t$. More generally, using \eqref{eq:walksinconeeq} and Corollary \ref{cor:albegraicity}, the generating function $\tilde{E}_{k,k_{1},k_{2}}(t)$ is algebraic when $3\nmid k_{1}-k_{2}$. When $3| k_{1}-k_{2}$ the generating function $E_{k,k_{1},k_{2}}(t)$ remains D-finite, however we can prove that it is not algebraic.

\subsection{Asymptotic number of walks in a cone}
Recall from Corollary \ref{cor:asymptotics} that we know the asymptotic form of the coefficients of $\tilde{E}(t,e^{i\alpha})$ for fixed $\alpha\in(0,\pi)$. By symmetry we have $\tilde{E}(t,e^{i\alpha})=\tilde{E}(t,e^{-i\alpha})$, so we can extend this result to $\alpha\in(-\pi,0)$. We also have similar results at $\alpha=0$ and $\alpha=\pi$. Using these results together with \eqref{eq:walksinconeeq} we derive the asymptotic form of the coefficients of $\tilde{E}_{k,k_{1},k_{2}}(t)$:
\[[t^n]\tilde{E}_{k,k_{1},k_{2}}(t)=-\frac{2\cdot 3^{5-\frac{6}{k_{2}-k_{1}}}\sin\left(\frac{k_{1}\pi}{k_{2}-k_{1}}\right)\sin\left(\frac{(k_{1}-2k)\pi}{k_{2}-k_{1}}\right)}{\pi(k_{2}-k_{1})^2\left(1+2\cos\left(\frac{2\pi}{k_{2}-k_{1}}\right)\right)\Gamma\left(-\frac{3}{k_{2}-k_{1}}\right)}n^{-1-\frac{3}{k_{2}-k_{1}}}3^{n}.\]
Writing $\beta=\frac{(k_{2}-k_{1})\pi}{3}$ to denote the opening angle of the cone, the polynomial term $n^{-1-\frac{\pi}{\beta}}$ above appears much more generally \cite{denisov2015random}, but to our knowledge the exact value of the constant term above is new.

\section{Conclusion and other lattices}
\label{sec:otherlattices}
In the previous sections we described how to count walks on the Kreweras lattice by length and winding number. As a consequence we counted walks on a conic subsection of this lattice using the reflection principle. The same method can be readily applied to extract analogous results for walks on the triangular lattice, the square lattice and the king lattice, shown in Figure \ref{fig:lattices}. On each lattice we count walks on cones with opening angle $Nr$ for any positive integer $N$, where $r=\frac{\pi}{3}$ for the triangular lattice and Kreweras lattice and $r=\frac{\pi}{4}$ for the square lattice and king lattice. On the square lattice, many of our results were found by Timothy Budd in 2017 using a very different method \cite{budd2017winding}.

One direction for further research is to generalise this method to apply to walks using different step sets $S\subset\{-1,0,1\}^{2}$, as have been studied on the quarter plane and, more recently, on the three quarter plane.

\acknowledgements{I would like to thank Mireille Bousquet-M\'elou for many helpful suggestions which greatly improved the manuscript.}

%% if you use biblatex then this generates the bibliography
%% if you use some other method then remove this and do it your own way
%\printbibliography

\end{document}